\documentclass[noid,thms,letterpaper,namedreferences]{kluwer}
\usepackage{amsfonts, amssymb, times, mathptm}
\def\copyrightline{}
\def\kaplogo{}

\makeatletter\@openrightfalse\def\@listI{%
  \leftmargin \leftmargini
  \topsep  9\p@ \@plus 3\p@ \@minus 5\p@
  \partopsep 3\p@ \@plus 1\p@ \@minus 2\p@
   \parsep 2\p@ \@plus 1\p@ \@minus 1\p@
\itemsep\parsep}
\@listI
\makeatother

\let\nmodels\nvDash

\let\copyrightline\relax
\let\kaplogo\relax
\let\impl\to
\def\fin{\mathrm{fin}}
\def\inf{\mathrm{inf}}

\def\T{\mathfrak{T}}
\def\Ht{\mathbf{Ht}}\def\H{\ensuremath{\mathbf{H}}}\def\L{\ensuremath{\mathbf{L}}}
\def\qp#1{\ensuremath{{#1}\mathbf{\pi{+}}}}
\def\Tree{\mathrm{Tree}}
\def\Prop{\mathrm{Prop}}
\def\Fin{\mathrm{Fin}}
\def\Arity{\mathrm{Arity}}

\def\SwS{\ensuremath{\mathrm{S}\omega\mathrm{S}}}
\def\S4{\ensuremath{\mathbf{S4}}}
\def\LC{\ensuremath{\mathbf{LC}}}
\def\G{\ensuremath{\mathbf{G}}}

\begin{document}\bibliographystyle{klunamed}
\begin{article}
\begin{opening}
\title{Decidability of Quantified Propositional\\
Intuitionistic Logic and S4 on Trees}
\author{Richard \surname{Zach}
 \email{rzach@ucalgary.ca}}
\institute{Department of Philosophy\\ University of Calgary\\
2500 University Drive NW\\
Calgary, Alberta T2N 1N4\\ Canada}
\date{Draft, March 17, 2002---Comments welcome!}

%\begin{ao}\\
%Richard Zach\\
%Department of Philosophy\\
%University of Calgary\\
%Calgary, Alberta  T2N 1N4\\
%Canada\\
%rzach@ucalgary.ca\\
%http://www.ucalgary.ca/$\sim$rzach/
%\end{ao}

\begin{abstract}
Quantified propositional intuitionistic logic is obtained from
propositional intuitionistic logic by adding quantifiers $\forall p$,
$\exists p$, where the propositional variables range over
upward-closed subsets of the set of worlds in a Kripke structure.  If
the permitted accessibility relations are arbitrary partial orders,
the resulting logic is known to be recursively isomorphic to full
second-order logic \cite{Kremer:97}.  It is shown that if the Kripke
structures are restricted to trees, the resulting logics are
decidable.  The result also transfers to modal \S4{} and some
G\"odel-Dummett logics with quantifiers over propositions.
\end{abstract}

\end{opening}

\section{Introduction}

Quantified propositional intuitionistic logic is obtained from
propositional intuitionistic logic by adding quantifiers $\forall p$,
$\exists p$ over propositions.  In the context of Kripke semantics, a
proposition is a subset of the worlds in a model structure which is
upward closed, i.e., if $h \in P$, then $h' \in P$ for all $h' \ge h$.
For propositional intuitionistic logic~\H, several classes of model
structures are known to be complete, in particular the class of all
partial orders, as well as the class of trees and some of its
subclasses.  When quantifiers over propositions are added, these
results no longer hold.  \inlinecite{Kremer:97} has shown that the
quantified propositional intuitionistic logic \qp{\H} based on the
class of all partial orders is recursively isomorphic to full
second-order logic.  He raised the question of whether the logic
resulting from restriction to trees is axiomatizable. The main part of
this note establishes that, in fact, it is decidable.

It should be pointed out right away that the trees we consider here
are all subtrees of the complete tree of height and arity~$\omega$.
That is, trees of uncountable arity, or height more than~$\omega$ are
excluded.  This is in accord with Kripke's \shortcite{Kripke:65}
intuitive interpretation of his possible world semantics for
intuitionistic logic.  In this interpretation, Kripke explains, the
worlds in a structure correspond to ``points in time (or `evidential
situations')'' and the accessibility relation $\le$ holds between
worlds $h$, $h'$ if ``as far as we know, at time $h$, we may later
gain enough information to advance to $h'$.''  If the language is
countable, then at each point, there are only countably many sentences
about which we could discover new information. So at each point, there
are only countably many possibilities for advancing to a new
evidentiary situation, i.e., the tree of evidentiary situations should
have arity $\le \omega$.  Allowing trees of transfinite height would
correspond, in this interpretation, to allowing a transfinite process
of gathering of evidence. A ``jump'' to a new evidentiary situation
only after an infinite amount of time and investigation seems counter
to the spirit of Kripke's interpretation; hence, trees should be of 
height at most~$\omega$.

The rest of this note is organized as follows: Section~2 introduces
the logics considered, and contains several observations regarding the
relationship between the classes of formulas valid on various classes
of trees.  Section~3 presents the decidability result for quantified
propositional intuitionistic logic.  Section~4 outlines how the
results transfer to a proof of decidability of modal \S4{} with
propositional quantification on similar types of Kripke structures.
(Propositionally quantified \S4{} on general partial orders is also
known to be not axiomatizable.)  Intermediate logics based on linear
orders (i.e., 1-ary trees), which correspond to G\"odel-Dummett
logics, are also considered.  A concluding section discusses
limitations and possible extensions of the method.

\section{Quantified propositional intuitionistic logics}

\begin{defn} An \emph{model structure}~$\langle g, K, \le\rangle$
is given by a set of worlds~$K$, an initial world $g \in K$, and a
partial order $\le$ on~$K$, for which $g$ is a least element.  Given
a structure, an \emph{(intuitionistic) proposition} is a subset
$P \subseteq K$ so that when $h \in P$ and $h' \ge h$, then also $h'
\in P$.  A \emph{valuation}~$\phi$ is a function mapping the
propositional variables to propositions of~$M$.  A \emph{model}~$M =
\langle g, K, \le, \phi\rangle$ is a structure together with a
valuation.  If $P$ is a proposition in the model $M$, then $M[P/p]$
is the model which is just like $M$ except that it assigns the
proposition $P$~to~$p$.
\end{defn} 

\begin{defn}
If $M = \langle g, K, \le, \phi\rangle$ is a model, $h \in K$, and $A$
is a formula, we define what it means for \emph{$A$ to be true at
$h$}, denoted $M, h \models A$, by induction on formulas as follows:
\begin{enumerate}
\item $M, h \models p$ if $h \in \phi(p)$; $M, h \nmodels \bot$.
\item $M, h \models B \land C$ if $M, h \models B$ and $M, h \models C$.
\item $M, h \models B \lor C$ if $M, h \models B$ or $M, h \models C$.
\item $M, h \models B \impl C$ if, for all $h' \ge h$, either $M, h'
\not\models B$ or $M, h' \models C$.
\item $M, h \models \forall p\, B$, if, for all
propositions $P$, $M[P/p], h \models B$.
\item $M, h \models \exists p\, B$ if there is a proposition~$P$ so that
$M[P/p], h \models B$.
\end{enumerate}
The constant $\bot$ is always assigned the empty proposition; $\neg A$
abbreviates $A \impl \bot$, hence, $M, h \models \neg B$ iff for all
$h' \ge h$, $M, h' \nmodels B$.
\end{defn} 

\begin{defn} Given a model~$M$ and a formula $A$, the \emph{proposition 
defined by $A$} is the set $M(A) = \{h : M, h \models A\}$.
\end{defn}

\begin{prop}
$M(A)$ is a proposition. In fact we have:
\[\begin{array}{rcl@{\qquad}rcl}
M(p) & = & \phi(p) &  M(\bot) & = & \emptyset\\ 
M(A \land B) & = &  M(A) \cap M(B) & 
M(A \lor B) & = & M(A) \cup M(B)\\ 
M(\forall p\, A) & = & \bigcap_P M[P/p] A &
M(\exists p\, A) & = & \bigcup_P M[P/p] A\\
\multicolumn{6}{c}{M(A \impl B)  =  \{h : \textrm{for\ all\ }h'
\ge h, \textrm{\ if\ }h' \in M(A) \textrm{\ then\ }h' \in M(B)\}}
\end{array}\]
\end{prop}

\begin{pf} By induction on the complexity of formulas.\qed
\end{pf}

\begin{defn} A model $M$ \emph{validates}~$A$, $M \models A$, if 
$M, g \models A$. A model structure $S$ validates $A$, if every model
based on $S$ validates $A$.  $A$ is \emph{valid in a class of model
structures~$\mathfrak{C}$}, $\mathfrak{C} \models A$, if $M \models A$
for all models~$M$ based on structures in~$\mathfrak{C}$. $A$ is
\emph{valid}, if $M \models A$ for any model~$M$.\end{defn}

\begin{defn} A tree~$T$ is a subset of $\omega^*$, the set of words 
over $\omega$, which is closed under initial segments.  $T$ is
partially ordered by the prefix ordering $\le$ defined as: $x \le y$
if $y = xz$ for some~$z$, and totally ordered by the lexicographic
order~$\preceq$.  The empty word~$\Lambda$ is the least element in both
orderings. The set $T_\omega = \omega^*$ itself is a tree, the
\emph{complete infinitary tree.}  The set $T_n = \{i: 0 \le i < n\}^*$
($n \le \omega$) is also a tree (called the \emph{complete $n$-ary
tree}). 
\end{defn}

\begin{defn}\label{logics}
 We consider the following classes of model structures on trees:
\[\begin{array}{rclrcl}
\T & = &\{\langle \Lambda, T, \le\rangle : T \textrm{ is a
tree}\} & 
\T_{n} & = & \{T_n\}, \\
\multicolumn{6}{c}{\T_\fin  = \{\langle \Lambda, T, \le\rangle : T \textrm{ is a finite tree}\}.}
\end{array}\]
These models tructures give rise to the following quantified propositional
logics:\[\begin{array}{rclrcl}
\qp{\H} & = & \{A : {}\models A\} &
\qp{\Ht} & = & \{A : \T \models A\} \\
\qp{\Ht_n} & = & \{A : \T_{n} \models A\} &
\qp{\Ht^\fin} & = & \{A : \T_\fin \models A\}. \\
\end{array}\]
\end{defn}

To each of these quantified propositional logics $\qp{\L}$ corresponds
a propositional logic~\L{} obtained by restriction to quantifier-free
formulas.  These all collapse to~$\H$, i.e., $\H = \Ht = \Ht_n
=\Ht^\fin,$ for $n \ge 2$ \cite{Gabbay:81}.  The \emph{quantified}
propositional logics, however, do not:

\begin{prop}\label{rels}
1. $\qp{\H} \subsetneq \qp{\Ht} \subsetneq \qp{\Ht_n}$ and $\qp{\Ht} \subsetneq \qp{\Ht^\fin}$.\\
2. $\qp{\Ht^\fin} \not\subseteq \qp{\Ht_n}$ and 
$\qp{\Ht_n} \not\subseteq \qp{\Ht^\fin}$.
\end{prop}

\begin{pf} The inclusions $\qp{\H} \subseteq \qp{\Ht} \subseteq \qp{\Ht_n}$, 
and $\qp{\Ht} \subseteq \qp{\Ht^\fin}$ are obvious.  

%To see $\qp{\Ht_\omega} \subseteq \qp{\Ht_n}$,
%consider a model $M = \langle \Lambda, T, \le,
%\phi\rangle$ with $T \in \T_n$.  Define a model $M' = \langle \Lambda,
%T_\omega, \le, \phi'\rangle$ as follows: $T_0 = T$ and $\phi_0 =
%\phi$.  To construct $T_{i+1}$, $\phi_{i+1}$ take all words $x$ of
%length~$i$.  $x$ has a greatest immediate successor $xn$ in the
%lexicographical ordering.  Let $T_{i+1} = T_i \cup \{ xmy : xny \in
%T_{i}, m > n\}$ and $\phi_{i+1}(p) = \phi_i(p) \cup \{xmy : xny \in
%\phi_i(p), m > n\}$.  Let $T' = \bigcup T_i$ and $\phi' = \bigcup
%\phi_i$.  It is easily seen that $M \models A[h]$ iff $M' \models
%A[h]$ for $h \in T$.  Thus, if $A \notin \qp{\Ht_n}$, then $A
%\notin\qp{\Ht_\omega}$.

To show that the first inclusion is proper, consider:
\begin{eqnarray*}
A &=& \forall p(\neg p \lor \neg\neg p) \impl 
\forall p\forall q((p \impl q)\lor(q\impl p))
\end{eqnarray*}
Then $\qp{\H} \nmodels A$: The 4-element diamond is a countermodel.
On the other hand, $\qp{\Ht} \models A$, since any $h$ with $h \models
\forall p(\neg p \lor \neg\neg p)$ is so that for all $h', h'' \ge h$,
either $h' \ge h''$ or $h'' \ge h'$.  To see this, suppose $h', h''
\ge h$ but neither $h' \le h''$ nor $h'' \le h'$.  Consider the
proposition $P = \{k: k \ge h'\}$.  Then $M[P/p], h' \models p$, and
hence $M[P/p], h \nmodels \neg p$.  On the other hand, $M[P/p], k
\nmodels p$ for any $k \ge h''$.  Hence, $M[P/p], h'' \models \neg p$
and so $M[P/p], h \nmodels \neg\neg p$.  In other words, the part of
the model above $h$ is linearly ordered, and so $h \models \forall
p\forall q((p \impl q) \lor (q \impl p))$.

For the second inclusion, take $B = \forall p(p \lor \neg p)$. Since
$\forall p(p \lor \neg p)$ is true at any $h$ which has no successor
worlds in a model (a leaf node) and false otherwise, $\neg B$ will be
true iff the model has no leaf node.  Since complete trees don't have
leaf nodes, $\qp{\Ht_n} \models \neg B$ but $\qp{\Ht} \nmodels
\neg B$.\footnote{This example is due to Tomasz Po\l{}acik.  Instead of $p
\lor \neg p$ one can use any classical tautology which is not
derivable in intuitionistic logic.}

On the other hand, in a finite tree, every branch has a world with no
successors. If $M$ is a model based on a finite tree, for every world
$h$ there is a world $h' \ge h$ such that $M, h' \models \forall p(p
\lor \neg p)$.  Hence, for every world $h$, $M, h \nmodels \neg B$ and
consequently $M, h \models \neg\neg B$.  Thus,
$$\begin{array}[b]{rclrcl}
\qp{\Ht^\fin} & \models & \neg\neg B,\textrm{ but} &
\qp{\Ht}, \qp{\Ht_n} & \nmodels & \neg\neg B; \\
\qp{\Ht^\fin} & \nmodels & \neg B,\textrm{ but} &
\qp{\Ht_n} & \models & \neg B.
\end{array}\eqno\Box$$
\end{pf}

\section{Decidability results}

\begin{thm}[\opencite{Kremer:97}] \qp{\H} is recursively isomorphic to full 
second-order logic. \end{thm}

\begin{thm} Each logic from Definition~\ref{logics}, except \qp{\H},
is decidable.
\end{thm}

\begin{pf}
We use Rabin's tree theorem \cite{Rabin:69}.  That theorem says that
\SwS, the monadic second-order theory of $T_\omega$, is decidable.
We reduce validity of quantified propositional formulas to truth of
formulas of \SwS.

The language of \SwS{} contains two relation symbols $\le$ and
$\preceq$, for the prefix ordering and the lexicographical ordering,
respectively, and a constant $\Lambda$ for the empty word.  
\emph{Finiteness} is definable in \SwS: $X$ is finite iff it
has a largest element in the lexicographic ordering~$\preceq$. Let $x
\le_1 y$ say that $y$ is an immediate successor of $x$. Then we have:
\begin{eqnarray*}
\Tree(T) & = & \Lambda \in T \land 
\forall x \in T\, \forall y (y \le x \impl y \in T)) \\
\Prop(T) & = & \forall x \in T\, \forall y (x \le y \impl y \in T)) \\
\Arity_n(T) & = & \forall x \in T\, \exists^{= n} y (x \le_1 y)) \textrm{ if $n < \omega$} \\
\Fin(T) & = & \exists x\forall y \in T\, y \preceq x)
\end{eqnarray*}
which say that $T$ is a tree (with root $\Lambda$), 
a proposition, has arity~$n$, and is finite,
respectively.

If $A$ is a formula of quantified propositional logic, define $A^x$
by:
\[\begin{array}{rclrcl}
p^x & = & x \in X_p & (B \impl C)^x & = & \forall y \in T(x \le y \impl (B^y \impl C^y)) \\
\bot^x & = & \bot & (\forall p\, B)^x & = & \forall X_p((X_p \subseteq T \land \Prop(X_p)) \impl B^x) \\
(B \land C)^x & = & B^x \land C^x & 
(\exists p\, B)^x & = & \exists X_p(X_p \subseteq T \land \Prop(X_p) \land B^x),\\
(B \lor C)^x & = & B^x \lor C^x
\end{array}\]
where $y$ is a new variable not previously used in the translation.
Now let
\begin{eqnarray*}
\Psi(A, \qp{\Ht}) & = & \forall T(\Tree(T) \impl A^x[\Lambda/x]) \\ 
\Psi(A, \qp{\Ht_n}) & = & \forall T((\Tree(T) \land
\Arity_n(T)) \impl A^x[\Lambda/x] \quad(n < \omega)\\ 
\Psi(A, \qp{\Ht^\fin}) & = & \forall T((\Tree(T) \land \Fin(T)) 
\impl A^x[\Lambda/x]) \\ 
\Psi(A, \qp{\Ht_\omega}) & = & \forall T(\forall z(z \in T) \impl 
A^x[\Lambda/x])
\end{eqnarray*}

We may assume, without loss of generality, that $A$ is closed (no free
propositional variables).

We have to show that $\SwS \models \Psi(A, \qp\L)$ iff $\qp\L \models
A$.  First, let $M = \langle \Lambda, K, \le, \phi\rangle$ be an
\qp\L-model (obviously, we may assume that $\Lambda$ is the root). If
$M, \Lambda \nmodels A$, then $M(A) \neq K$.  Define a variable
assignment~$s$ for second-order variables by $s(T) = K$.  Then it is
easy to see that $M(A) = \{x \in K: \SwS \models A^x [s]\}$.  Thus,
$\Psi(A, \qp\L)$ is false in~\SwS.

Conversely, if $\SwS \nmodels \Psi(A, \qp\L)$, then there is a
counterexample witness~$X$ for the initial universal quantifier
$\forall T$, which is a tree (in the respective class), $\Lambda \in
X$, and $\SwS \nmodels A^x[\Lambda/x] [s]$ for $s(T) = X$.  (For the
case of $\L = \qp{\Ht_\omega}$, $X = T_\omega$.)

We show that for any $s$ with $s(T) = X$, the model $M = \langle
\Lambda, X, \le, \phi\rangle$ with $\phi(p) = s(X_p)$ is such
that $M(A) = \{x \in X: \SwS \models A^x [s]\}$.  This is obvious if $A
= p$, $A = B \land C$ or $A = B \lor C$.  Suppose $A = B \impl C$.
Then $x \in M(A)$ iff for all $y \in X$ with $x \le y$, $y \notin
M(B)$ or $y \in M(C)$.  By induction hypothesis, $y \notin M(B)$ iff
$\SwS \nmodels B^y [s]$; similarly for $y \in M(C)$. So $x \in M(A)$
iff $\SwS \models A^x [s]$.  If $A = \forall p\, B$, then $x \in M(A)$
iff for all propositions~$P$ in $X$, $x \in M[P/p](B)$.  This is the
case, by induction hypothesis, iff for all upward-closed subsets~$P$
of $X$, $\SwS \models B^x [s']$ where $s'$ is like $s$ except $s'(X_p)
= P$; but this is true just in case $\SwS \models \forall X_p((X_p
\subseteq T \land \Prop(X_p)) \impl B^x)$. (Similarly for the case of
$A = \exists p\,B$.) Hence, if $A$ is closed and $\SwS \nmodels
\Psi(A, \qp\L)$, the structure $M = \langle \Lambda, X, \le,
\phi\rangle$ is a countermodel for $A$.\qed
\end{pf}

\section{S4 and G\"odel-Dummett logics}

Modal logic \S4{} is closely related to intuitionistic logic, and its
Kripke semantics is likewise based on partially ordered structures and
trees.  In the modal context, a proposition is any (not necessarily
upward-closed) subset of the set of worlds.  Adding quantifiers over
propositions to \S4, we obtain the logic \qp{\S4}.  

Specifically, the semantics of \qp{\S4} is like that for \qp{\H},
except that an \emph{\S4-proposition} in $M$ is a subset $P \subseteq
K$, and valuations $\phi$ map variables to \S4-propositions.  We have
the two modal operators $\Box$ and $\Diamond$.  $M, h \models A$ 
is then defined by
\begin{enumerate}
\item $M, h \models p$ if $h \in \phi(p)$; $M, h \nmodels \bot$.
\item $M, h \models B \land C$ if $M, h \models B$ and $M, h \models C$.
\item $M, h \models B \lor C$ if $M, h \models B$ or $M, h \models C$.
\item $M, h \models B \impl C$ if $M, h \nmodels B$ or $M, h \models C$.
\item $M, h \models \Box B$ if all $h' \ge h$, $M, h' \models B$.
\item $M, h \models \Diamond B$ if some $h' \ge h$, $M, h' \models B$.
\item $h \models \forall p\, B$, if, for all
propositions $P$, $M[P/p], h \models B$.
\item $h \models \exists p\, B$ if there is a proposition~$P$ so that
$M[P/p], h \models B$.
\end{enumerate}
Depending on the class of Kripke structures considered, we obtain
logics \qp{\S4}, \qp{\mathbf{S4t}}, \qp{\mathbf{S4t}_n},
\qp{\mathbf{S4t}^\fin} (for the class of partial orders, trees,
$n$-ary trees, and finite trees, respectively).

The McKinsey-Tarski $T$-embedding of \H{} into \S4 \cite[Theorem
5.1]{McKinseyTarski:48} can be straightforwardly extended to the propositional quantifiers. For a formula $A$ in the language of \qp{\H},
define a formula $A^T$ of \qp{\S4} as follows:
\[
\begin{array}{rcl@{\qquad}rcl}
p^T & = & \Box p & (B \impl C)^T  & = & \Box(B^T \impl C^T)\\
\bot^T & = & \Box \bot & (\forall p\, B)^T & = & \forall p\, B^T \\
(B \land C)^T & = & B^T \land C^T & (\exists p\, B)^T & = & \exists p\, B^T\\
(B \lor C)^T & = & B^T \lor C^T
\end{array}\]

\begin{prop}
$\qp{\H} \models A$ iff $\qp{\S4} \models A^T$.
\end{prop}

\begin{pf}
Let $M = \langle g, K, \le, \phi\rangle$ be an intuitionistic
structure, and suppose $M, h \nmodels A$.  Consider the \S4-structure
$M' = \langle g, K, \le, \phi'\rangle$ with $\phi'(p) = \phi(p)$.  By
induction on the complexity of formulas, $M', h \nmodels A^T$.

Conversely, if $M' = \langle g, K, \le, \phi'\rangle$ is an
\S4-structure and $M', h \nmodels A^T$, then $M'', h \nmodels A^T$,
where $M'' = \langle g, K, \le, \phi''\rangle$ with $\phi''(p) =
M'(\Box p)$. \qed
\end{pf}

Note that the order structure of $M$ and $M'$ was not changed, so the
result holds also relative to any class of tree structures.  We can
therefore obtain separation results like those in
Proposition~\ref{rels} for the propositionally quantified variants of
\S4{} by considering the images under the $T$-embedding of the
formulas $A$, $\neg B$, and $\neg\neg B$ from the proof of
Proposition~\ref{rels}.

\inlinecite{Fine:70} and \inlinecite{Kremer:93} showed that \qp{\S4},
like \qp{\H} is not axiomatizable.  By the same method used above, the
decidability of \qp{\S4} can be established if one is only interested
in trees.

\begin{prop}
 \qp{\mathbf{S4t}}, \qp{\mathbf{S4t}_n}, and \qp{\mathbf{S4t}^\fin}
are decidable.
\end{prop}

\begin{pf}
We change the definition of $A^x$ as follows:
\[\begin{array}{rclrcl}
p^x & = & x \in X_p & (\Diamond B)^x & = & \exists y(x \le y \land B^y) \\
\bot^x & = & \bot & (\Box B)^x & = & \forall y(x \le y \impl B^y) \\
(B \land C)^x & = & B^x \land C^x & \forall p\,B^x & = & \forall X_p(X_p \subseteq T  
\impl B^x) \\
 (B \lor C)^x & = & B^x \lor C^x & \exists p\,B^x & = & \exists X_p(X_p \subseteq T \land B^x) \\
(B \impl C)^x & = & B^x \impl C^x
\end{array}\]
(where $y$ is new.) The definition of $\Psi(A, \qp\L)$ and the proof
that $\SwS \models \Psi(A, \qp\L)$ iff $\qp\L \models A$ ($\L$ one of
\qp{\mathbf{S4t}}, \qp{\mathbf{S4t}_n}, \qp{\mathbf{S4t}^\fin}) is the
same as for the intuitionistic case, mutatis mutandis.\qed
\end{pf}

Other logics which can be treated using the method used above are
G\"odel-Dummett logics. These logics were originally characterized as
many-valued logics over subsets of~$[0, 1]$. Here, a \emph{valuation}
is a mapping of propositional variables to truth values. A
valuation~$v$ is extended to formulas by:
\[
\begin{array}{cc}
\begin{array}{rcl}
v(\bot)                         &=& 0 \\
v(A \land B)     &=& \min(v(A), v(B))
\end{array} &
\begin{array}{rcl}
v(A \lor B)      &=& \max(v(A), v(B))\\
v(A \impl B)  &=&
\left\{ \begin{array}{cl}
   1            & {\rm if\ } v(A) \leq v(B) \\
   v(B)         & {\rm otherwise}
  \end{array}\right.
\end{array}
\end{array}
\]
In the quantifier-free case, taking any infinite subset of $[0, 1]$ as
the set of truth values results in the same set of tautologies,
axiomatized by $\LC = \H + (A \impl B) \lor (B \impl A)$.  This is no
longer the case if we add propositional quantifiers.  In the
many-valued context, these can be introduced by:
\begin{eqnarray*}
v(\exists p\, A) &=& \sup \{ v[w/p](A) : w \in V \}\\
v(\forall p\, A) &=& \inf \{ v[w/p](A) : w \in V\},
\end{eqnarray*}
where $v[w/p]$ is the valuation which is like $v$ except that it
assigns the value~$w$ to $p$.  The resulting class of tautologies
depends on the order structure of $V \subseteq [0, 1]$. In fact, there
are $2^{\aleph_0}$ different propositionally quantified
G\"odel-Dummett logics. %\cite{BaazVeith:98}.

\LC{} is also characterized as the set of formulas valid on the
infinite 1-ary tree $\T_1$.  The G\"odel-Dummett logic which
corresponds to this characterization is $\G_\downarrow\pi$ based on
the truth-value set $V_\downarrow = \{0\} \cup \{1/n : n \ge 1\}$,
i.e., $\G_\downarrow\pi = \qp{\Ht_1}$ (\opencite{BaazZach:98},
Proposition 2.8). The intersection of all finite-valued
G\"odel-Dummett logics, however, coincides with $\G_\uparrow\pi$ with
truth value set $V_\uparrow = \{1\} \cup \{1-1/n : n \ge 1\}$.  Since
$\G_\uparrow\pi \neq \G_\downarrow\pi$, this shows that the formulas
valid on the infinite 1-ary tree is not identical to the class of
formulas valid on all 1-ary trees of finite height.  This latter logic
was studied and axiomatized by \inlinecite{BCZ:00}.

\section{Conclusion}

As noted in the introduction, the notion of trees we consider is
the only one which accords with Kripke's intuitive interpretation of
intuitionistic model structures.  It might nevertheless be interesting
to consider more general classes of trees (i.e., partial orders with
least element and where $h \not\le h'$ and $h' \not\le h$ guarantees
that for no $g$ is $h, h' \le g$), or well-founded trees (every branch
is well-ordered).  

The problem of the complexity of the resulting quantified
propositional logics on such structures, however, remains open.  It is
not known whether the monadic second-order theory of such partial
orders is decidable, in fact, it most likely is not.  If it were,
however, the reduction given here would immediately yield the
decidability results for the quantified propositional logics on such
structures.

We can also easily obtain further decidability results for logics
based on classes of trees which are definable in the language of \SwS.
This includes, e.g., tress of finite arity, trees of finite height, and
trees of arity or height $\le n$ for some~$n$.

\section*{Acknowledgements}
Thanks to Chris Ferm\"uller, Yuri Gurevich, Tomasz Po{\l}acick, Saharon
Shelah, and Wolfgang Thomas for comments and helpful discussion.

\end{article} \end{document}